\documentclass{amsart}

\newtheorem{theorem}{Theorem}[section]
\newtheorem{lemma}[theorem]{Lemma}

\theoremstyle{definition}

\theoremstyle{remark}

\numberwithin{equation}{section}

\begin{document}

\title{\textbf{Classifying convex compact ancient solutions to the affine curve shortening flow}}

\author{Shibing Chen}
\address{Department of Mathematics, University of Toronto, Toronto, Ontario
Canada M5S 2E4}
\curraddr{Department of Mathematics, University of Toronto, Toronto, Ontario
Canada M5S 2E4}
\email{sbchen@math.toronto.edu}
\thanks{}


\date{February 10}

\dedicatory{}


 \begin{abstract}
 In this paper we classify convex compact ancient solutions to the affine curve shortening flow: namely, any convex compact ancient solution to the affine curve shortening flow must be a shrinking ellipse. The method combines a rescaling argument inspired by \cite{Wang}, affine invariance of the equation and monotonicity of the affine isoperimetric ratio. We will give two proofs.  The essential ideas are related, but
the first one uses level set represenation
of the evolution.
 The second proof employs Schauder's estimates, and it also provides a new simple proof for the corresponding classification result to the higher dimensional affine normal flow.
 \end{abstract}
 
 \maketitle

 \section{Introduction}
The affine curve shortening flow (ACSF) is described by a time-dependent embedding of a fixed smooth closed curve $\ell$: $X(\cdot,t): \ell\rightarrow \mathbb{R}^{2}$ satisfying the following
evolution equation
  \begin{equation}
\label{e2}
\frac{\partial X}{\partial t}=\kappa^{\frac{1}{3}}\textbf{n},
\end{equation}
where $\textbf{n}$ is the inner normal of the evolving curve, and $\kappa$ is the curvature the evolving curve.
ACSF was studied by many authors. In particular Sapiro and Tannenbaum \cite{ST} proved that when the initial curve is convex, it will converge to an elliptic point under the affine curve shortening flow. Later Andrews \cite{An1} generalized their result to the higher dimensional case, which is called the affine normal flow. One of the most interesting results for ACSF was proved by Angenent, Sapiro and Tannenbaum
\cite{AST}, where they showed that any simple closed curve will shrink to an elliptic point under ACSF. For all these results one can find
corresponding results in mean curvature flow, see \cite{GH}, \cite{Huisken} and \cite{Grayson}. ACSF and its higher dimensional analog have very interesting  applications in geometric analysis, for example it has been used by Andrews \cite{An1} to give a very beautiful proof of the affine isoperimetric inequality. Besides its mathematical interest, ACSF is one of the fundamental equations in image processing (see \cite{AGLM}), and it has been extensively used in this subject, we refer the reader to \cite{FC} and references therein.

This paper concerns the classification of convex compact ancient solutions of ACSF. An ancient solution to ACSF is a solution with existence time $(-\infty, T).$ For the curve shortening flow, Wang \cite{Wang} got the first partial classification result of the convex ancient solution, he proved that if a convex ancient solution to the curve shortening flow sweeps the whole space $\mathbb{R}^{2}$, it must be a shrinking circle, otherwise the convex ancient solution must be defined in a strip region and he indeed constructed such solutions by some compactness argument. Later, Daskalopoulos, Hamilton and Sesum \cite{DHS} showed that apart from the shrinking circle, the so called ¡±Angenent oval¡± (a convex ancient solution of the curve shortening flow discovered by Angenent that decomposes into two translating solutions of the flow) is the only other embedded convex compact ancient solution of the curve shortening flow.
For the higher dimensional affine normal flow, the classification of convex compact ancient solution was done by Loftin and Tsui in \cite{LT}, where they proved the only compact, ancient solutions to the affine normal flow in $\mathbb{R}^{n}$, $n\geq 3$, are ellipsoids. However, their method can not be used for the ACSF, since it is based on showing the cubic tensor vanishes which only implies the hypersurface is a hyperquadric in higher
dimensions. In this paper, by combining a rescaling argument inspired by \cite{Wang}, affine invariance of the equation and monotonicity of the affine isoperimetric ratio, we classify convex compact ancient solutions of ACSF, namely a convex compact ancient solution to ACSF must be a shrinking ellipse.

We represent the affine curve shortening flow using the level set equation:
\begin{equation}\label{e1}
\text{div}(\frac{Du}{|Du|})=\frac{1}{|Du|^{3}}.
\end{equation}
Let $u:\mathbb{R}^{2}\rightarrow \mathbb{R}$ be a solution of \eqref{e1}, then the level set $\{u=-t\}$ evolves under ACSF.
 By the affine invariance of ACSF, we have that for any unimodular affine transformation $A$, $u(Ax)$ is still
an solution of \eqref{e1}. A convex compact ancient solution to the affine curve shortening flow corresponds to a complete convex solution to
\eqref{e1} with compact level set (see Lemma \ref{l1} and the discussion after it). In the following we use the rescaling $u_{h}(x)=\frac{1}{h}u(h^{\frac{3}{4}}x),$ and it is easy to see that $u_{h}$ is still an solution of $\eqref{e1}$. We will also assume $u\geq 0$, and $u(0,0)=0$.

\begin{theorem}\label{t1}
A convex compact ancient solution to the affine curve shortening flow must be a shrinking ellipse. Equivalently, a complete convex solution to
\eqref{e1} with compact level set must be $\frac{3}{4}|x|^{\frac{4}{3}}$ up to an affine transformation.
\end{theorem}

\section{First proof of the main theorem}
First, we will prove a lemma which is used to show that once we represent an ancient convex solution to \eqref{e2} as a solution $u$ to the level set
equation \eqref{e1}, $u$ must be a convex function.
\begin{lemma}\label{l1}
Let $\Omega$ be a smooth, bounded, convex domain in $\mathbb{R}^{2}$. Let $u$ be the solution of \eqref{e1}, vanishing on $\partial\Omega$.
Then for any constant $h$ satisfying $\inf_{\Omega}u<h<0$, the level set $\{u=h\}$ is convex. Moreover, $\log(-u)$ is a concave function.
\end{lemma}
\emph{Proof.} Observe $\varphi:=-\log(-u)$ satisfies
$$|D\varphi|^{2}\displaystyle{\sum_{i,j=1}^{2}}(\delta_{ij}-\frac{\varphi_{i}\varphi_{j}}{|D\varphi|^{2}})\varphi_{ij}=e^{\frac{1}{3}\varphi}.$$
Since $\varphi(x)\rightarrow +\infty$ as $x\rightarrow\partial\Omega$, the result in \cite{K}(Theorem 3.13) implies $\varphi$ is convex.

With the above lemma and the Lemma 4.4 in \cite{Wang}, we know that any convex compact ancient solution to the affine curve shortening flow can be represented as a convex solution $u$ to equation  \eqref{e1}.

\emph{Proof of Theorem \ref{t1}.}
Step 1. We will find a sequence of levels $h_{m}\rightarrow\infty$ as $m\rightarrow\infty$, such that $\{u=h_{m}\}$ looks more and more like an ellipse in the sense that $(1-\epsilon_{m})E_{m}\subset \{u=h_{m}\}\subset (1+\epsilon_{m})E_{m}$ (here ``$\subset$" means ``enclosed by", and the dilation $(1+\epsilon_{m})E_{m}$ is always done with respect to the center of the ellipse ) with the constants $\epsilon_{m}\rightarrow 0$, where $E_{m}$ is a sequence of ellipses.

First we take any sequence $h_{m}\rightarrow\infty$ as $m\rightarrow\infty$. By John's lemma we can find an unimodular affine transformation $A_{m}$, so that $\{x|u(A_{m}x)=h_{m}\}$ has a good shape in the sense that $B_{m}\subset\{x|u(A_{m}x)=h_{m}\}\subset 2B_{m}$, where $B_{m}$ is a ball and here the subscript $m$
does not mean the radius. Then we blow down the solution $u(A_{m}x)$ to $v_{m}(x)=\frac{1}{h_{m}}u(A_{m}(h_{m}^{\frac{3}{4}}x))$, which is still a solution of \eqref{e1}.

Now $\{v_{m}=1\}$ has good shape. If we interpret the level curves of $v_{m}$ evolve under the affine curve shortening flow, it will take time 1 for it to shrink to the origin, which combines with the good shape property of $\{v_{m}=1\}$ imply that the width of $\{v_{m}=1\}$ has upper bound independent of $m$.
By Blaschke's selection theorem, passing to a subsequence we can assume $\{v_{m}=1\}$ converges to a limiting bounded convex curve $\ell_{0}$. So the affine curve shortening flow with intial data $\{v_{m}=1\}$  converges to the flow  $\ell_{t}$ ($t$ is from 0 to 1) with initial data $\ell_{0}$ in the hausdorff distance.
Since affine curve shortening flow shrinks a convex compact curve to an elliptic point, the curve $\ell_{t}$ looks more and more like an ellipse as $t\rightarrow 1$ in the sense that $(1-\epsilon_{t})E_{t}\subset \ell_{t} \subset (1+\epsilon_{t})E_{t}$ (here and below, all ellipses centered at the origin) with $\epsilon_{t}$ going to 0, where $E_{t}$ is a family of ellipses. First, we fix a sequence $a_{i}\rightarrow 0$ as $i\rightarrow \infty$.
By the convergence of the flows with intial data $\{v_{m}=1\}$  to the flow $\ell_{t}$ with initial data $\ell_{0}$, for every $a_{i}$ we can choose $m_{i}$ large enough to guarantee that
$\{v_{m_{i}}=a_{i}\}$ is very close to $\ell_{1-a_{i}}$ in hausdorff distance (close enough so that $(1-2\epsilon_{1-a_{i}})E_{1-a_{i}}\subset \{v_{m_{i}}=a_{i}\} \subset (1+2\epsilon_{1-a_{i}})E_{1-a_{i}}$
and $a_{i}h_{m_{i}}\rightarrow\infty$ as $i\rightarrow\infty$. Now, rescaling back and after an affine transformation we see that $\{u=a_{i}h_{m_{i}}\}$ looks more and more like an ellipse as $i\rightarrow\infty$ in the sense as we discussed at the beginning of this step.

Step 2. Now we use the sequence $h_{m}$ we found at step 1 to blow down and renormalize the solution.
Let $A_{m}$ be the unimodular affine transformation so that $A_{m}^{-1}$ turns $E_{m}$ to a circle. Let $v_{m}=\frac{1}{h_{m}}u(A_{m}(h_{m}^{\frac{3}{4}}x))$.
One can see that $\{v_{m}=1\}$ converges to a unit circle $S^{1}$ centered at the origin in hausdorff distance. So the convex functions $v_{m}$ converges
to $\frac{3}{4}|x|^{\frac{4}{3}}$ locally uniformly in the unit disk $D_{1}$. In particular the gradient of $v_{m}$ converges to the gradient of $\frac{3}{4}|x|^{\frac{4}{3}}$ uniformly in the Disk $D_{\frac{3}{4}}$ (here $\frac{3}{4}$ is the radius of the disk) as $m\rightarrow\infty$, which means the value of $|Dv_{m}|$ on $\{v_{m}=\frac{1}{2}\}$ converges uniformly to a constant (the norm of the gradient of $\frac{3}{4}|x|^{\frac{4}{3}}$ at $\{\frac{3}{4}|x|^{\frac{4}{3}}=\frac{1}{2}\}$) as $m\rightarrow\infty$. Then by the equation \eqref{e1}, we see the curvature of $\{v_{m}=\frac{1}{2}\}$ (given by the left hand side of \eqref{e1}) converges uniformly to a constant (the curvature of $\{\frac{3}{4}|x|^{\frac{4}{3}}=\frac{1}{2}\}$) as $m\rightarrow\infty$.

Recall that the affine isoperimetric ratio of a convex curve $\mathcal{C}$
(parametrized by gauss map) is defined by $V^{-\frac{1}{3}}\int_{S^{1}}\kappa^{-\frac{2}{3}}d\theta$, where $V$ is the volume enclosed by the curve, and $\kappa$ is the curvature. Andrews \cite{An1} proved that the affine isoperimetric ratio is increasing to $2\pi^{\frac{2}{3}}$ for affine curve shortening flow with convex compact initial data. Now by the above uniform convergence of the curvature of $\{v_{m}=\frac{1}{2}\}$ to the curvature of the circle $\{\frac{3}{4}|x|^{\frac{4}{3}}=\frac{1}{2}\}$ and the hausdorff convergence of $\{v_{m}=\frac{1}{2}\}$ to the circle $\{\frac{3}{4}|x|^{\frac{4}{3}}=\frac{1}{2}\}$, we see that the affine isoperimetric ratio of $\{v_{m}=\frac{1}{2}\}$ converges uniformaly $2\pi^{\frac{2}{3}}$ as $m\rightarrow\infty$. Since $h_{m}\rightarrow\infty$, after blow-down and renormalization, $\{u=1\}$ becomes a level curve of $v_{m}$ below $\{v_{m}=\frac{1}{2}\}$, so the affine isoperimetric ratio of $\{u=1\}$ is greater than the affine isoperimetric ratio of $\{v_{m}=\frac{1}{2}\}$, which converges to $2\pi^{\frac{2}{3}}$, but by the affine isoperimetric inequality we know the affine isoperimetric ratio of $\{u=1\}$ is less or equal to $2\pi^{\frac{2}{3}}$, which means the affine isoperimetric ratio of $\{u=1\}$ is $2\pi^{\frac{2}{3}}$, namely it is an ellipse.
It is easy to see that we can apply the above argument for any level $\{u=h\}$, namely any level set of $u$ must be an ellipse.

\section{Second proof of the main theorem}
This proof includes the case for higher dimensional affine normal flow. Recall that the affine normal flow is described by a time-dependent embedding of a fixed $n$-dimensional hypersurface $M$, namely, $X(\cdot,t): M\rightarrow \mathbb{R}^{n+1}$ satisfying the following
evolution equation
  \begin{equation}
\label{e3}
\frac{\partial X}{\partial t}=K^{\frac{1}{n+2}}\textbf{n},
\end{equation}
where $\textbf{n}$ is the inner normal of $X(M,t)$, and $K$ is the Gauss curvature of the evolving hypersurface.  Note that \eqref{e3} has rescaling invariance, namely if
$X(\cdot,t)$ is a solution of \eqref{e3}, then $\lambda^{-\frac{n+2}{2n+2}}X(\cdot, \lambda t)$ ($\lambda$ is any positive constant) is also a solution. The equation is also invariant under unimodular affine transformation, namely if $X(\cdot,t)$ is a solution of \eqref{e3}, then $AX(\cdot,t)$ is also a solution for any unimodular affine transformation $A$.

This paper is concerned only with ancient convex compact solutions to \eqref{e3}, namely a solution with existence time starting from $-\infty$, ending at $0$. Below we will assume the solution shrinks to an elliptic point at origin when time goes to $0$.
 First, we take a sequence of positive numbers $T_{i}\rightarrow\infty $ as $i\rightarrow\infty$. Then by John's lemma, for every $i$, one can find a sequence of unimodular affine transformations $A_{i}$, so that $B_{i}\subset A_{i}(X(M\times\{-T_{i}\}))\subset (n+1)B_{i}$, where $B_{i}$ is a ball ($i$ is only an index, but not the radius). By the invariance properties discussed above, we see that $\widetilde{M}_{i}^{t}:=T_{i}^{-\frac{n+2}{2n+2}}A_{i}(X(M\times \{T_{i}t\}))$  evolves under the affine normal flow starting at time $-1$, ending at time $0$, with initial data $\widetilde{M}_{i}^{-1}=T_{i}^{-\frac{n+2}{2n+2}}A_{i}(X(M\times \{-T_{i}\}))$ which has good shape by the choice of $A_{i}$. By using $\widetilde{B}_{i}:=T_{i}^{-\frac{n+2}{2n+2}}B_{i}$ as a barrier, the good shape property combines with the assumption that it takes time $1$ for $\widetilde{M}_{i}^{-1}$ shrinking to origin force the radius of $\widetilde{B}_{i}$ to be bounded above by some constant independent of $i$.
By blaschke's selection theorem, after passing to a subsequence we can assume $\widetilde{M}_{i}^{-1}\rightarrow\widetilde{M}$ for some convex hypersurface $\widetilde{M}$. Denote $\widetilde{M}^{t}$ as a solution to the affine normal flow with initial data $\widetilde{M}$.

Now, for fixed $t\in(-1,0)$, by applying Schauder estimates for solutions of uniformly parabolic equations concave in the second derivatives as in \cite{An2} Section 8.1, we have that $\widetilde{M}_{i}^{t}\rightarrow \widetilde{M}^{t}$ in $C^{\infty}$ norm. Recall the affine isoperimetric ratio of a convex hpersurface $N$ (parametrized by gauss map) is given by $I(N):=V^{-\frac{n}{n+2}}\int_{S^{n}}K^{-\frac{n+1}{n+2}}d\Omega,$ where $V$ is the volume enclosed by $N$.
 The main property for affine isoperimetric ratio we will use below is that it is increasing under the affine normal flow (see \cite{An1}). So $I(\widetilde{M}_{i}^{t})\rightarrow I(\widetilde{M}^{t})$ as $i\rightarrow\infty.$ Since $-\frac{1}{T_{i}}>t$ when $i$ is large, we have
$$I(X(M\times\{-1\}))=I(\widetilde{M}_{i}^{-\frac{1}{T_{i}}})\geq I(\widetilde{M}_{i}^{t})\rightarrow I(\widetilde{M}^{t})\ \text{as}\ i\rightarrow \infty,$$
where the $``\geq"$ is because of the monotonicity of affine isoperimetric ratio under ACSF.

Hence by affine isoperimetric inequality, $$I(\text{ellipsoid})\geq I(X(M\times\{-1\}))\geq I(\widetilde{M}^{t})\rightarrow I(\text{ellipsoid})\ \text{as}\ t \rightarrow 0,$$ which implies
$I(X(M\times\{-1\}))=I(\text{ellipsoid}),$ Since equality is achieved in the affine
isoperimetric inequality only by ellipsoids, $X(M\times\{-1\})$ must be an ellipsoid. It is easy to see that the above argument works for $X(M\times\{-h\})$ with positive $h$, namely $X(M\times\{-h\})$ must be an ellipsoid for any positive $h$.

\emph{Remark.} It is easy to see that this method also works for $p$-centro affine normal flows, which was recently studied by Ivaki and Stancu \cite{IS}.
The result is `` The only compact, origin-symmetric, ancient solutions to the p-flow are shrinking ellipsoids."
Note that one can apply schauder estimates similarly as above, and should replace the affine isoperimetric ratio with p-affine isoperimetric ratio.
\bibliographystyle{amsplain}

\end{document}